\documentclass[envcountsect,12pt]{llncs}
\usepackage{fancyheadings}
\usepackage{amssymb,bbm,showlabels}
\usepackage{hyperref}
\usepackage[numeric]{amsrefs}
\makeatletter

\renewcommand{\BibLabel}{%
     \hyper@anchorstart{cite.\CurrentBib}\relax\thebib.\hyper@anchorend%
}
\makeatother

\usepackage{amsmath}
\usepackage{amsfonts}
\ifx\pdftexversion\undefined
  \usepackage[dvips]{graphicx}
\else
  \usepackage[pdftex]{graphicx}
\fi

\newcommand{\thmref}[1]{Theorem~\ref{#1}}
\newcommand{\lemref}[1]{Lemma~\ref{#1}}
\newcommand{\secref}[1]{Section~\ref{#1}}

\numberwithin{equation}{section}
         
\newcommand\g{\gamma}

\newcommand\e{\varepsilon}
\renewcommand\l{\lambda}

\newcommand\G{\Gamma}
\newcommand\f{\frac}

\newcommand{\Z}{{\mathbb{Z}}}
\newcommand{\R}{{\mathbb{R}}}

\newcommand{\C}{{\mathbb{C}}}

\renewcommand\i{^{-1}}
\renewcommand\({\left(}
\renewcommand\){\right)}

\newcommand{\ignore}[1]{}
\newcommand{\myignore}[1]{}
\newcommand{\mymyignore}[1]{}
\newcommand{\bx}{\hfill$\square$\vspace{.6cm}}
\begin{document}

\title{Spectral Analysis of Pollard Rho Collisions}

\date{April 17, 2006}
\author{Stephen D. Miller\thanks{Partially supported by NSF grant
DMS-0301172 and an Alfred P. Sloan Foundation Fellowship.
}\inst{1}~\,and\, Ramarathnam Venkatesan\inst{2}}

\institute{Einstein Institute of Mathematics \\ The Hebrew
University \\ Givat Ram, Jerusalem 91904, Israel \\ and \\
Department of Mathematics \\ Rutgers University \\ Piscataway, NJ
08854, USA \\
\email{miller@math.huji.ac.il}   \vspace{.2cm} \and   Microsoft Research \\
Cryptography and Anti-piracy Group \\ 1 Microsoft Way, Redmond, WA
98052, USA  \\ and \\
Cryptography Research Group \\ Microsoft Research India \\
Scientia - 196/36 2nd Main \\
Sadashivnagar, Bangalore 560 080, India     \\
 \email{venkie@microsoft.com}}

\maketitle

\begin{abstract} We show that  the classical Pollard $\rho$
algorithm for discrete logarithms produces a collision in expected
time $O(\sqrt{n}(\log n)^3)$. This is the first nontrivial
rigorous estimate for the collision probability for the {\em
unaltered} Pollard $\rho$ graph, and is close to the conjectured
optimal bound of $O(\sqrt{n})$. The result is derived by showing
that the mixing time for the random walk on this graph is $O((\log
n)^3)$; without the squaring step in the Pollard $\rho$ algorithm,
the mixing time would be exponential in $\log n$.  The technique
involves a spectral analysis of \emph{directed} graphs, which
captures the effect of the squaring step. \ignore{
 This is
necessary  because the well known equivalence of eigenvalue
separation and rapid mixing for {\em undirected} graphs
unfortunately breaks down for directed graphs.}

\vspace{.2 cm}

 Keywords:
Pollard Rho algorithm, discrete logarithm,
 random walk, expander graph, collision time, mixing time, spectral
analysis.
\end{abstract}

\section{Introduction}

Given a finite cyclic group $G$ of order $n$ and a generator $g$,
the Discrete Logarithm Problem (\textsc{dlog}) asks to invert the
map $y\mapsto g^y$ from $\Z/n\Z$ to $G$.  Its presumed difficulty
serves as the basis for several cryptosystems, most notably the
Diffie-Hellman key exchange and some elliptic curve cryptosystems.
Up to constant factors, the Pollard $\rho$ algorithm is the most
efficient and the only version with small memory known for solving
\textsc{dlog} on  a general cyclic group -- in particular for the
group of points of an elliptic curve over a finite field.

We quickly recall the algorithm now. First one randomly partitions
$G$ into three sets $S_1$, $S_2$, and   $S_3$.  Set $x_0=h$, or more
generally to a random power $g^{r_1} h^{r_2}$.   Given $x_k$, let
$x_{k+1}=f(x_k)$, where $f:G\rightarrow
  G$ is defined by
\begin{equation}\label{iteration}
    f(x) \ \ = \ \ \left\{
                     \begin{array}{ll}
                       gx\,, &  \ x\,\in\,S_1\,; \\
                       hx\,, & \ x\,\in\,S_2\,; \\
                       x^2\,, & \ x\,\in\,S_3\,.
                     \end{array}
                   \right.
\end{equation}
Repeat until a collision of values of the $\{x_k\}$ is detected
(this is done using Floyd's method of comparing $x_k$ to $x_{2k}$,
which has the advantage of requiring minimal storage).  We call
the underlying directed graph in the above algorithm (whose
vertices are the elements of $G$, and whose edges connect each
vertex $x$ to $gx$, $hx$, and $x^2$) as the {\em Pollard $\rho$
Graph}. At each stage $x_k$ may be written as $g^{a_ky+b_k}$,
where $h=g^y$. The equality of $x_k$ and $x_\ell$ means $a_{k\,} y
+ b_k = a_{\ell\,} y + b_\ell$, and solving for $y$ (if possible)
recovers the \textsc{dlog} of $h=g^y$.

The above algorithm {\em heuristically} mimics a random walk. Were
that indeed the case, a collision would be found in time
$O(\sqrt{n})$, where $n$ is the order of the group $G$.  (The
actual constant is more subtle; indeed, Teske \cite{teske1} has
given evidence that the walk is somewhat worse than random.)

The main result of this paper is the first rigorous nontrivial upper
bound on the collision time.  It is slightly worse than the
conjectured $O(\sqrt{n})$, in that its runtime is
$\widetilde{O}(\sqrt{n})$, i.e.~off from $O(\sqrt{n})$ by at most a
polynomial factor in $\log n$.  As is standard and without any loss
of generality, we tacitly make the following
\begin{equation}\label{primeassumption}
    \textbf{assumption:~} \text{~the order~~}|G|\,=\,n\text{~~is
    prime.}
\end{equation}

\begin{theorem}\label{mainthm}  Fix $\e>0$.  Then the Pollard $\rho$
algorithm for discrete logarithms on $G$ finds a collision in time
$O_\e(\sqrt{n}\,(\log n)^3)$ with probability at least $1-\e$, where
the probability is taken over all partitions of $G$ into three sets
$S_1$, $S_2$, and $S_3$.
\end{theorem}

 In the black-box group model (i.e.~one
which does not exploit any special properties of the encoding of
group elements), a theorem of Shoup \cite{shoup} states that any
\textsc{dlog} algorithm needs $\Omega(\sqrt{n})$ steps. Hence, aside
from the probabilistic nature of the above algorithm and the extra
factor of $(\log n)^3$, the estimate of  Theorem~\ref{mainthm} is
sharp.

It should be noted that finding a collision does not necessarily
imply finding a solution to \textsc{dlog}; one must also show the
resulting linear equation is nondegenerate.  Since $n=|G|$ is prime
this is believed to happen with overwhelming probability, much more
so than for the above task of finding a collision in $O(\sqrt{n})$
time. This was shown for a variant of the Pollard $\rho$ algorithm
in \cite{horwen}, but the method there does not apply to the
original
 algorithm itself.  Using more refined techniques we are able
 to analyze this question
further; the results of these investigations will be reported upon
elsewhere.

 This paper is the first analysis of the unmodified
Pollard $\rho$ Graph, including the fact that it is
 {\em directed}.
    One can obtain the required rapid mixing result for directed
    graphs
by (a) assuming that  rapid mixing holds for the undirected
version, and
 (b) adding self-loops to each vertex.
 However, one still needs to prove (a), which
 in our situation is no simpler.  In addition,
 the loops and loss of direction
  cause short cycles, which lead to awkward
 complications in the context of studying collisions.

Technically, analyzing directed graphs from a spectral point of
view has the well known difficulty that a spectral gap is not
equivalent to rapid mixing.  A natural generalization of the
spectral gap is the operator norm gap of the adjacency matrix,
which suffices for our purposes (see \secref{rmdir}).  For a
recent survey of mixing times on directed graphs, see
\cite{montteta}.

The Pollard $\rho$ graph is very similar to the graphs introduced by
the authors in \cite{mv3}. These graphs, which are  related to
expander graphs, also connect group elements $x$ to $f(x)$ via the
operations given in (\ref{iteration}) -- in particular they combine
the operations of multiplication and squaring. The key estimate, a
spectral bound on the adjacency operator on this graph, is used to
show its random walks are rapidly mixing.  Though the Pollard $\rho$
walk is only {\em pseudorandom} (i.e., $x_{k+1}$ is determined
completely from $x_k$ by its membership in $S_1$, $S_2$, or $S_3$),
we are solely interested here in proving that it has a collision.
The notions of random walk and pseudorandom walk (with random
assignments of vertices in the sets $S_i$) coincide until a
collision occurs.

\subsection{Earlier Works}

Previous experimental and theoretical studies of the Pollard $\rho$
algorithm and its generalizations all came to the (unproven)
conclusion that it runs in $O(\sqrt{n})$ time; this is in fact the
basis for estimating the relative bit-for-bit security of elliptic
curve cryptosystems compared to others, e.g.~RSA.
 For an analysis of \textsc{dlog}
algorithms we refer the reader to the survey by Teske \cite{teske2},
and for an analysis of random walks on abelian groups, to the one by
Hildebrand \cite{hildebrand}.  For the related Pollard $\rho$
algorithm for factoring integers, Bach \cite{bach} improved the
trivial bound of $O(n)$ by logarithmic factors.

An important statistic of the involved graphs is the {\em mixing
time} $\tau$, which loosely speaking is the amount of time needed
for the random walk to converge to the uniform distribution, when
started at an arbitrary node.\footnote{There are many inequivalent
notions of mixing time (see \cite{lovasz}). Mixing time is only
mentioned for purposes of rough comparison between different
graphs; whatever we need about it is proved directly. Similarly,
the reader need not recall any facts about expander graphs, which
are mentioned only for motivation.} The existing approaches to
modeling Pollard $\rho$ can be grouped into two categories:

\begin{enumerate}

\item \emph{Birthday attack in a totally random model:}  each step
is viewed as a move to a random group element, i.e.~a completely
random walk. In particular one assumes that the underlying graph
has mixing time $\tau=1$ and that its degree equals the group
size; in reality the actual Pollard $\rho$ graph has degree only
3.  The $O(\sqrt{n})$ collision time is immediate for random walks
of this sort.

\item \emph{Random walk in an augmented graph}: The Pollard $\rho$
graph is modified by increasing the number of generators $k$, but
removing the squaring step. One then models the above transitions
as random walks on directed abelian Cayley graphs. To ensure the
mixing time is $\tau=O(\log |G|)$, however, the graph degree must
grow at least logarithmically in $|G|$.  The importance of $\tau$
stems from the fact that, typically, one incurs a overhead of
multiplicative factor of $\tau^{const}$ in the overall algorithm.

\end{enumerate}

Teske \cite{teske1}, based on Hildebrand's results
\cite{hildebrand} on random walks on the cyclic group $\Z/m\Z$
with respect to steps of the form $x\mapsto x+a_i$, $i\le k$,
shows that the mixing time of an algorithm of
 the second type
 is
on the order of  $n^{\frac{2}{k-1}}$; she gives supporting
numerics of random behavior for $k$ large. In particular, without
the squaring step the Pollard $\rho$ walk would have mixing time
on the order of  $n^2$, well beyond the expected $O(\sqrt{n})$
collision time.  This operation is an intriguing and
cryptographically\footnote{In this version one can derive a secure
hash function \cite{horwitz}  whose security is based on the
difficulty of the discrete logarithm problem; here the input
describes the path taken in the graph from a fixed node, and the
hash value is the end point.} important aspect of the
 Pollard $\rho$ algorithm, and makes it
 inherently {\it non-abelian}: the Pollard $\rho$ graphs
are not isomorphic to any abelian Cayley graphs.  Its effect
cannot be accounted for by any analysis which studies only the
additive structure of $\Z/m\Z$.

The present paper indeed analyzes the {\em exact} underlying
Pollard $\rho$ graph, without any modifications.  We are able to
show that the inclusion of the squaring step reduces the
 mixing
time $\tau$ from exponential in $\log n$, to $O((\log n)^3)$ --- see
the remark following Proposition~\ref{prhoseparation}.

 Our result and
technique below easily generalize from the unmodified Pollard $\rho$
algorithm, which has only 2 non-squaring operations, to the
generalized algorithms proposed by Teske \cite{teske1} which involve
adding further such operations.  Furthermore, it also applies more
generally to additional powers other than squares. We omit the
details, since the case of interest is in fact the most difficult,
but have included a sketch of the argument at the end of the paper.

\section{Rapid mixing on directed graphs}
\label{rmdir}

In the next two sections we will describe some results in graph
theory which are needed for the proof of \thmref{mainthm}. Some of
this material is analogous to known results for {\em undirected}
graphs (see, for example, \cite{bollobas}); however, since  the
literature on spectral analytic aspects of directed graphs is
relatively scarce, we have decided to give full proofs for
completeness.

The three properties of subset expansion, spectral gap, and rapid
mixing are all equivalent for families of undirected graphs with
fixed degree.   This equivalence, however, fails for directed
graphs. Although a result of Fill \cite{Fill} allows one to deduce
rapid mixing on directed graphs from undirected analogs, it involves
adding self-loops (which the Pollard $\rho$ graph does not have) and
some additional overhead.  In any event, it requires proving an
estimate about the spectrum of the undirected graph.   We are able
to use the inequality \cite[(A.10)]{mv3}, which came up in studying
related undirected graphs, in order to give a bound on the operator
norm of the directed graphs.  This bound, combined with
Lemma~\ref{mixlem}, gives an estimate of $\tau = O((\log n)^3)$ for
the mixing time of the Pollard $\rho$ graph.

Let $\G$ denote a graph with a finite set of vertices $V$ and edges
$E$.  Our graphs will be directed graphs, meaning that each edge has
an orientation;  an edge from $v_1$ to $v_2$ will be denoted by
$v_1\rightarrow v_2$. Assume that $\G$ has {\em degree} $k$, in
other words that each vertex has exactly $k$ edges coming in and $k$
edges coming out of it. The {\em adjacency} operator $A$ acts on
$L^2(V)=\{f:V\rightarrow\C\}$ by summing over these $k$ neighbors:
\begin{equation}\label{adjacency}
    (Af)(v) \ \ = \ \ \sum_{v\rightarrow w} \,f(w)\,.
\end{equation}
Clearly constant functions, such as $\mathbbm{1}(v)\equiv 1$, are
eigenfunctions of $A$ with eigenvalue $k$.  Accordingly,
$\mathbbm{1}$ is termed the trivial eigenfunction and $k$ the
trivial eigenvalue of $A$.  Representing $A$ as a $|V|\times|V|$
matrix, we see it has exactly $k$ ones in each row and column,
with all  other entries equal to  zero. It follows that
$\mathbbm{1}$ is also an eigenfunction with eigenvalue $k$ of the
adjoint operator $A^*$
\begin{equation}\label{adjoint}
    (A^*f)(v) \ \ = \ \ \sum_{w\rightarrow v} \,f(w) \, ,
\end{equation} and that all eigenvalues $\l$ of $A$ or $A^*$
satisfy the bound $|\l|\le k$.

The subject of {\em expander graphs} is concerned with bounding the
 (undirected) adjacency operator's restriction to the
subspace $L_0=\{f \in L^2(V) \mid f \perp {\mathbbm 1}\}$,
i.e.~the orthogonal complement of the constant functions under the
$L^2$-inner product. This is customarily done by bounding the
nontrivial eigenvalues away from $k$. However, since the adjacency
operator $A$ of a directed graph might not be self-adjoint, the
operator norm can sometimes be a more useful quantity to study. We
next state a lemma relating it to the rapid mixing of the random
walk. To put the statement into perspective, consider the $k^r$
random walks on $\G$ of length $r$ starting from any fixed vertex.
One expects a uniformly distributed walk to land in any fixed
subset $S$ with probability roughly $\f{|S|}{|V|}$. The lemma
gives a condition on the operator norm for this probability to in
fact lie between $\f 12 \f{|S|}{|V|}$ and $\f32 \f{|S|}{|V|}$ for
moderately large values of $r$.  This can alternatively be thought
of as giving an upper bound on the mixing time.

\begin{lemma}\label{mixlem}
Let $\G$ denote a directed graph of degree $k$ on $n$ vertices.
Suppose that there exists a constant $\mu<k$ such that $\| A f \|
\le \mu\| f \|$ for all $f\in L^2(V)$ such that $f\perp
\mathbbm{1}$. Let $S$ be an arbitrary subset of $V$.  Then the
number of paths of length $r \ge \f{\log(2n)}{\log(k/\mu)}$ which
start from any given vertex and end in $S$ is between $\f12 k^r
\f{|S|}{|V|}$ and $\f32 k^r \f{|S|}{|V|}$.
\end{lemma}
\begin{proof}
Let $y$ denote an arbitrary vertex in $V$, and $\chi_S$ and
$\chi_{\{y\}}$ the characteristic functions of $S$ and $\{y\}$,
respectively.  The  number of paths of length $r$ starting at $y$
and ending in $S$ is exactly the $L^2(V)$-inner product $ \langle
\chi_S , A^r \chi_{\{y\}}\rangle$.  Write
\begin{equation}\label{orthsplit}
    \chi_S  \ \ = \ \ \f{|S|}{n}\mathbbm{1} \, + \, w \ \ \ \
    \text{and} \ \ \ \ \chi_{\{y\}} \ \ = \ \ \f{1}{n}\mathbbm{1} \,
    + \, u\, ,
\end{equation}
where $w,u\perp \mathbbm{1}$.  Because $\mathbbm{1}$ is an
eigenfunction of $A^*$, $A$ preserves the orthogonal complement of
$\mathbbm{1}$, and thus
\begin{equation}\label{AwAubdd}
    \|A^r u \| \ \ \le \ \ \mu\,\|A^{r-1}u\| \ \ \le \ \ \cdots \ \
    \le
    \ \ \mu^{r\,}\|u\|\,.
\end{equation}
Also, by orthogonality
\begin{equation}\label{uwnormbd}
    \| w \| \  \le  \ \| \chi_S \| \ = \ \sqrt{|S|} \ \ \ \ \
    \text{and} \ \ \ \ \  \|u\| \  \le \ \| \chi_{\{y\}} \| \
= \ 1\,.
\end{equation}
We have that $A^r \chi_{\{y\}} = \f{1}{n}k^r\mathbbm{1} + A^r
u$, so the inner product may be calculated as
\begin{equation}\label{innerprodcalc}
   \langle
\chi_S , A^r \chi_{\{y\}}\rangle  = \ \ \f{|S|}{n}k^r \, + \,
\langle w,A^ru\rangle\,.
\end{equation}
It now suffices to show that the absolute value of the second term
on the righthand side is bounded by half of the first term.  Indeed,
\begin{equation}\label{errbd}
|\langle w,A^ru\rangle| \ \ \le \ \ \|w\|\,\|A^r u\| \ \ \le \ \
\mu^{r}\sqrt{|S|}\,,
\end{equation}
 and
\begin{equation}\label{ineqchain}
    \mu^{r}\sqrt{|S|}  \  \le
    \  \f{1}{2n}  k^r
\sqrt{|S|}  \ \le  \    \f{1}{2}k^r\f{|S|}{n}
\end{equation}
when $r \ge \f{\log(2n)}{\log(k/\mu)}$. \bx
\end{proof}

\section{Collisions on the Pollard $\rho$ graph}

In this section, we prove an operator norm bound on the Pollard
$\rho$ graph that is later used in conjunction with
\lemref{mixlem}. These graphs are closely related to an undirected
graph studied in \cite[Theorem 4.1]{mv3}. We will start by quoting
a special case of the key estimate of that paper, which concerns
quadratic forms.    At first glance, the analysis is reminiscent
of the of the Hilbert inequality from analytic number theory (see
\cite{montgomery,steele}), but where the quadratic form
coefficients are expressed as $1/\sin(\mu_j-\mu_k)$.

Let $n$ be an odd integer and $\l_k = |\cos(\pi k/n)|$ for $k\in
\Z/n\Z$.  Consider the quadratic form $Q:\R^{n-1}\rightarrow \R$
given by
\begin{equation}\label{quadform}
    Q(x_1,\ldots,x_{n-1}) \ \ := \ \
    \sum_{k\,=\,1}^{n-1}\,x_k\,x_{2k}\,\l_{k}\,,
\end{equation}
in which the subscripts are interpreted modulo $n$.
\begin{proposition}\label{qformbd}
There exists an absolute constant $c>0$ such that
\begin{equation}\label{qformbdineq}
    |Q(x_1,\ldots,x_{n-1})| \ \ \le \ \ \( 1-\f{c}{(\log
    n)^2} \) \,\sum_{k\,=\,1}^{n-1}x_k^2 \, .
\end{equation}
\end{proposition}
\begin{proof}
Let $\g_k$ be arbitrary positive quantities (which will be
specified later in the proof).  Since
\begin{equation}\label{gammaineq}
 \g_k \,x_k^2 \ + \  \g_k\i\,x_{2k}^2 \  \pm \
 2 \,x_k \, x_{2k} \ \  = \  \   \(\g_k^{1/2}\,x_k \
 \pm \ \g_k^{-1/2}\,x_{2k}\)^2 \\  \ge \ \ 0 \, ,
\end{equation}
one has that
\begin{equation}\label{boundq}
    |Q({\vec x})| \  \ \le  \  \
\f{1}{2}\,\sum_{k\,=\,1}^{n-1}\(\g_k\,x_k^2\ + \
\g_k\i\,x_{2k}^2\)\,\l_k \ \ =   \ \
\f{1}{2}\,\sum_{k\,=\,1}^{n-1}x_k^2\(\g_k \,\l_k + \g_{\bar{2}k}\i
\l_{\bar{2}k}\),
\end{equation}
where $\bar{2}$ denotes the multiplicative inverse to 2 modulo $n$.
The proposition follows if we can choose $\g_k$ and an absolute
constant $c>0$ such that
\begin{equation}\label{critbd}
\g_k\,\l_k \ + \ \g_{\bar{2}k}\i \, \l_{\bar{2}k} \ \ < \ \  2 \ -
\ \f{c}{(\log n)^2} \ \ \  \ \ \text{ for all~}\,1 \,  \le \,  k
\, < \, n\,.
\end{equation}

Now we come to the definition of the $\g_k$. We set $\g_k=1$ for
$n/4 \le k \le 3n/4$; the definition for the set of other nonzero
indices $\cal S$ is more involved. For $s\ge 0$, define
$$t_s \ = \ 1 \, - \, s \, \f{d}{(\log n)^2}\,,$$ where $d>0$
is a small constant that is chosen at the end of the proof. Given an
integer $\ell$ in the range $-n/4 < \ell < n/4$, we define $u(\ell)$
to be order to which 2 divides $\ell$. For the       residues
$k\in\cal S$, which are all equivalent modulo $n$ to some integer
$\ell$ in the interval $-n/4 < \ell < n/4$, we define
$\g_k=t_{u(\ell)}$. Note also that $\l_k \le 1/\sqrt{2}$ for
$k\notin {\cal S}$, and is always $\le 1$. With these choices the
lefthand side of (\ref{critbd}) is bounded by
\begin{equation}\label{4choices}
    \g_k\,\l_k \ + \ \g_{\bar{2}k}\i \, \l_{\bar{2}k} \ \ \le \ \
\left\{
  \begin{array}{ll}
    \f{1}{\sqrt{2}}+ \f{1}{\sqrt{2}}\, , &  \ \ k, \bar{2}k\notin {\cal S} \\
    \f{1}{\sqrt{2}} + \g_{\bar{2}k}\i , &  \ \ k \notin {\cal S},\,  \bar{2}k\in {\cal S} \\
   \g_k +  \f{1}{\sqrt{2}} \, , &
 \ \ k\in {\cal S}, \, \bar{2}k\notin{\cal S}  \\
   \g_k + \g_{\bar{2}k}\i   , & \ \ k,\bar{2}k\in {\cal S}.
  \end{array}
\right.
\end{equation}
In the last case, the residues $k$ and $\bar{2}k$ both lie in
$\cal S$.  The integer $\ell \equiv \bar{2}k\! \pmod n$, $-n/4 <
\ell < n/4$, of course satisfies the congruence $2\ell\equiv
k\!\pmod n$. Since $k\in {\cal S}$,  $2\ell$ is the unique integer
in $(-n/4,n/4)$ congruent to $k$.  That means $\g_{k}=t_{s+1}$ and
$\g_{\bar{2}k}=t_{s}$ for some positive integer $s=O(\log n)$. A
bound for the last case in (\ref{4choices}) is therefore
$t_{s+1}+t_s\i=2-d/(\log n)^2 + O(s^2d^2/(\log n)^4)$.  We
conclude in each of the four cases that, for $d$ sufficiently
small, there exists a positive constant $c>0$ such that
(\ref{critbd}) holds.
 \bx
\end{proof}

The Pollard $\rho$ graph, introduced earlier, is the graph on
$\Z/n\Z$ whose edges represent the possibilities involved in
applying the iterating function (\ref{iteration}):
\begin{equation}\label{prhograph}
\gathered \text{$\G$ has vertices $V=\Z/n\Z$ and directed edges
$x\rightarrow x+1$, $x\rightarrow x+y$, }
\\
\text{and $x\rightarrow 2x$ for each  $x\in V$  (where $y\neq 1$) .}
\endgathered
\end{equation}

\begin{proposition}\label{prhoseparation}
Let $A$ denote the adjacency operator of the graph (\ref{prhograph})
and assume that $n$ is prime. Then there exists an absolute constant
$c>0$ such that
\begin{equation}\label{prhobd}
    \| Af \| \ \ \le \ \ \( 3 - \f{c}{(\log n)^2}\) \, \| f\|
\end{equation}
for all $f\in L^2(V)$ such that $f\perp {\mathbbm 1}$.
\end{proposition}
\begin{proof}
Let $\chi_{k}:\Z/n\Z\rightarrow\C$ denote the additive character
given by $\chi_k(x)=e^{2\pi i k x/n}$.  These characters, for $1 \le
k< n$, form a basis of functions $L_0=\{f\in L^2\mid f\perp
{\mathbbm 1}\}$. The action of $A$ on this basis is given by
\begin{equation}\label{Aaction}
    A\,\chi_k \  \ = \
      \ d_k \, \chi_k  \, + \,  \chi_{2k} \ \  \ , \ \ \ \ \text{where}
    \ \ d_k \ = \  e^{2\pi i k/n} \, + \, e^{2\pi i k y/n}\,.
\end{equation}
One has that $|d_{k}|\,=\,2|\cos(\f{\pi k (y-1)}{n})|=2\l_{k(y-1)}$.
Using the inner product relation
\begin{equation}\label{innerprodofchar}
    \langle \chi_{k},\chi_{\ell} \rangle  \ \ = \ \ \left\{
\begin{array}{ll}
 n\,, & \ \
k=\ell \\
   0\ , & \ \ \hbox{otherwise\,,}
     \end{array}
      \right.
\end{equation}
we compute that $\|f\|^2=n\sum |c_k|^2$,  where $f=\sum_{k\neq
0}c_k\chi_k$.  Likewise,
\begin{multline}\label{normcomputation} \|Af\|^2 \  \ = \
\ \langle Af , Af \rangle  \ \ = \\ \sum_{k,\ell\neq 0}
c_k\,\overline{c_{\ell}}\,\left[\langle d_k
\chi_k,d_\ell\chi_\ell\rangle + \langle\chi_{2k},\chi_{2\ell}\rangle
+ \langle d_k \chi_k,\chi_{2\ell}\rangle+ \langle
\chi_{2k},d_\ell\chi_\ell
\rangle\right] \\
\le  \ \ n\(5\sum |c_k|^2 + 2\sum |c_k||c_{2k}||d_{2k}|\).
\end{multline}
Note that $|d_k|=2\l_{k(y-1)}$,  and  that $y-1$ and 2 are
invertible in $\Z/n\Z$, by assumption in (\ref{prhograph}).  The
result now follows from (\ref{qformbdineq}) with the choice of
$x_{2(y-1)k}=|c_k|$.

 \bx
\end{proof}

\noindent {\bf Remark:} the above Proposition, in combination with
Lemma~\ref{mixlem}, is the source of the $\tau = O((\log n)^3)$
mixing time estimate for the Pollard $\rho$ graph that we
mentioned in the introduction.

\begin{proof}[of Theorem~\ref{mainthm}]
Consider the set $S$ of the first $t=\lfloor \sqrt{n} \rfloor$
iterates $x_1,x_2,\ldots,x_{t}$.  We may assume that $|S|=t$, for
otherwise a collision has already occurred in the first $\sqrt{n}$
steps. Lemma~\ref{mixlem} and Proposition~\ref{prhoseparation}
show that the probability of a walk of length $r\gg (\log n)^3$
reaching $S$ from any fixed vertex is at least $1/(2\sqrt{n})$.
Thus the probabilities that
$x_{t+r},x_{t+2r},x_{t+3r},\ldots,x_{t+kr}$ lie in $S$ are all,
independently, at least $1/(3t)$.  One concludes that for $k$ on
the order of $3bt$,  $b$ fixed, the probability that none of these
points lies in $S$ is at most $(1-\f{1}{3t})^{3bt} \approx
e^{-b}$, which is less than $\e$ for large values of $b$.
\end{proof}

\noindent {\bf Generalizations}: the analysis presented here extends
to generalized Pollard $\rho$ graphs in which each vertex $x$ is
connected to others of the form $xg_i$, for various group elements
$g_i$, along with powers $x^{r_j}$.  This can be done as follows.
First of all, if $r$-th powers are to be used instead of squares,
then the subscript $2k$ in (\ref{quadform}) must be changed to $rk$.
The key bound on (\ref{qformbdineq}), stated here for $r=2$, in fact
holds for any fixed integer $r>1$ which is relatively prime to $n$
\cite[Appendix]{mv3}. Thus changing the squaring step to
$x\rightarrow x^r$  does not change the end results. Secondly, the
proof of the bound (\ref{prhobd}) requires only some cancellation in
(\ref{normcomputation}).  If additional operations are added, the
cross terms from which the cancellation was derived here are still
present.  Thus Proposition~\ref{prhoseparation} is remains valid,
only with the 3 replaced by the degree of the graph.  Provided this
degree (= the total number of operations) is fixed, the graph still
has rapid mixing.

It is unclear if including  extra power operations speeds up the
discrete logarithm algorithm.  However, the rapid mixing of such
random walks may have additional applications, such as to the
stream ciphers in \cite{mv3}.

\vspace{.2cm}
 \noindent
 {\bf Acknowledgements:} the authors wish to thank R.
 Balasubramanian,
Michael Ben-Or, Noam Elkies, David Jao, L\' aszl\' o Lov\' asz,
and Prasad Tetali for helpful discussions and comments.

\end{document}